\documentclass[12pt]{amsart}

\usepackage{amsmath,amssymb,enumerate}
\usepackage{relsize}%

\usepackage{epsfig,fancyhdr,color}

\usepackage{amssymb}
\usepackage{amsmath,amsthm}
\usepackage{bbm}
\usepackage{latexsym}
\usepackage{amscd}
\usepackage{psfrag}
\usepackage{graphicx}
\usepackage[all]{xy}

\usepackage{verbatim}
\usepackage{amssymb}
\usepackage{amsmath}
\usepackage{graphicx}
\usepackage{color}
\usepackage{amsthm}
\usepackage{tikz}
\usetikzlibrary{arrows,positioning,decorations.pathmorphing, decorations.markings} 
\def\e{\epsilon}

\usepackage{rotating}
\usepackage{adjustbox}

\usepackage{xcolor}

\usepackage[font=small,labelfont=bf]{caption}

\usepackage{mathtools}

\usepackage{tikz-cd}

\usepackage{stackrel}

\usepackage[utf8]{inputenc}
\usepackage{amsmath}
\usepackage{amssymb}

\usepackage{mathabx,epsfig}
\def\acts{\mathrel{\reflectbox{$\righttoleftarrow$}}}

\usepackage[utf8]{inputenc}

\numberwithin{equation}{section}

\definecolor{NoteColor}{rgb}{1,0,0}

\newtheorem{theorem}{\rm\bf Theorem}[section]

\newtheorem{corollary}[theorem]{\rm\bf Corollary}
\newtheorem*{theorem 1}{\rm\bf Proposition 1}
\newtheorem*{theorem 2}{\rm\bf Proposition 2}

\theoremstyle{definition}

\theoremstyle{remark}

\def\interieur#1{\mathord{\mathop{\kern 0pt #1}\limits^\circ}}

\def\interieur#1{\mathord{\mathop{\kern 0pt #1}\limits^\circ}}

\def\hyperp{{\rm I}\kern-.3ex{\rm H}}

\usepackage{mathabx,epsfig}
\def\acts{\mathrel{\reflectbox{$\righttoleftarrow$}}}

\newtheorem*{theorem*}{Theorem}

\begin{document}

\title[The action P(SL(2,${\mathbb Z}$)) $\acts ~{\mathcal Tess}^+$ ]{Universal Spin Teichm\"uller Theory, I. \\The action of P(SL(2,${\mathbb Z}$)) on ${\mathcal Tess}^+$}
\
\author{Robert Penner} \address {\hskip -2.5ex Institut des Hautes \'Etudes Scientifiques\\ 35 route des Chartres\\ Le Bois Marie\\ 91440 Bures-sur-Yvette\\ France\\ {\rm and}~Mathematics Department, UCLA\\ Los Angeles, CA 90095\\USA} \email{rpenner{\char'100}ihes.fr}

\thanks{It is a pleasure to thank Igor Frenkel, Athanase Papadopoulos and Anton Zeitlin for discussions
and Barry Mazur and Dennis Sullivan for remarks on exposition.
}

 \date{\today}


\begin{abstract} 
Earlier work took as universal mapping class group the collection
${\rm PPSL}(2,{\mathbb Z})$
of all piecewise ${\rm PSL}(2,{\mathbb Z})$ homeomorphisms of the unit circle $S^1=\partial{\mathbb D}$ with finitely many breakpoints among the rational points in $S^1$.  
The spin mapping class group P(SL(2,${\mathbb Z}$))
introduced here consists of all piecewise-constant maps $S^1\to{\rm SL}(2,{\mathbb Z})$ which projectivize to an element of ${\rm PPSL}(2,{\mathbb Z})$.   We also introduce a spin universal Teichm\"uller space ${\mathcal Tess}^+$ covering the earlier 
universal Teichm\"uller space ${\mathcal Tess}$ of tesselations of ${\mathbb D}$ with fiber the space
of ${\mathbb Z}$/2 connections on the graph dual to the tesselation in ${\mathbb D}$.
There is a natural action ${\rm P(SL(2},{\mathbb Z}))\acts{\mathcal Tess}^+$ which is universal for finite-type hyperbolic surfaces with spin structure in the same sense that ${\rm PPSL}(2,{\mathbb Z})\acts{\mathcal Tess}$
is universal for finite-type hyperbolic surfaces.  
Three explicit elements of P(SL(2,${\mathbb Z}$))
are defined combinatorially via their actions on
${\mathcal Tess}^+$, and the main new result here is 
that they generate P(SL(2,${\mathbb Z}$)).
Background, including
material on hyperbolic and spin structures on finite-type surfaces, is sketched down to first principles in order to 
motivate the new constructions and to provide an overall survey.  A companion paper to this one gives a finite presentation
of the universal spin mapping class group ${\rm P(SL(2},{\mathbb Z}))$ introduced here.

\end{abstract}
\maketitle






\setcounter{footnote}{0} 

\section*{Introduction}

I am happy and honored to contribute here to the 90th birthday Festschrift for Valentin  Po\'enaru, that is, 
for Po.    We made our acquaintance only relatively recently with my residence in Bures-sur-Yvette
and quickly became fast friends. We have enjoyed wide-ranging discussions across math and physics,
in both scientific and sociological detail, as well as our fair share of delectations of fine wines and dinners including
Milen and their lovely family.  Po and Milen are a magical couple of transcendent love, their lives
comfortably merging Po's science with Milen's graphic art, even as this confluence
 has been passed on to the next generation of their children.

Po is one of the main founders of the Orsay topology group,  and thanks to  his network of friends and relations, he is  the one in that group who was the most open to the world mathematical community, both East and West.
He has always worked on the most difficult problems \cite{3+4}, including the Poincar\'e Conjecture in dimension three
and the smooth Schoenflies Problem in dimension four, on which his studies continue today, including the related
massive project  \cite{QSF1,QSF2,QSF3} in geometric group theory.

Along with his lifelong best friend Barry Mazur,  independently
Po \cite{pocell} and Barry \cite{barry} gave the headline first examples of pseudo 4-cells which are not
topological 4-cells, but whose product with the unit interval are
topological 5-cells. 
(A pseudo cell is a contractible compact combinatorial manifold with boundary.)
Among other fruitful collaborations
including his celebrated Immersion Theorem \cite{Po+} with Andre Haefliger,
Po's work with Boone and Haken \cite{Po++}
provides a fundamental tool for decision problems in topology,

Meanwhile, Po the man is the gentlest, kindest and warmest person I know, a beacon of calm and clear reason
to whom I often turn for guidance both mathematical and otherwise.   Here to Po in gratitude, admiration and friendship
in this celebratory volume, I offer the universal action ${\rm P(SL(2},{\mathbb Z}))\acts {\mathcal Tess}^+$ of the spin 
mapping class group on its spin Teichm\"uller space as a healtfelt 90th birthday present.
In order to formulate this, I first take two paragraphs to introduce the 
main new characters and recall their antecedents:  

\medskip

\noindent {\bf Universal Spaces.} First define a {\it tesselation} $\tau$ of the Poincar\'e disk
${\mathbb D}$ to be a locally finite collection of geodesics decomposing ${\mathbb D}$ into ideal triangles.  We shall typically choose a distinguished oriented edge or simply {\it doe} in the tesselation $\tau$, and define the auxiliary space
$${\mathcal Tess}'=\{{\rm tesselations~of~}{\mathbb D}~{\rm with~doe}\},$$
which we shall prove later is homeomorphic to the space of orientation-preserving homeomorphisms of the circle
$S^1=\partial {\mathbb D}$.  The M\"obius group ${\rm PSL}(2,{\mathbb R})$ acts on ${\mathcal Tess}'$ with Fr\'echet-manifold quotient
$${\mathcal Tess}={\mathcal Tess}'/{\rm PSL}(2,{\mathbb R})$$
the {\it universal Teichm\"uller space} of \cite{unicon}.  The spin version is introduced here as the
total space of the bundle $${\mathcal Tess}^+\to{\mathcal Tess}$$ of all equivalence classes
of finitely supported functions $\tau\to \{0,1\}$, where the equivalence relation is generated by
changing all three values on the boundary of a fixed triangle complementary to $\tau$ in ${\mathbb D}$.

\medskip

\noindent {\bf Universal Groups.}
Recall that PSL(2,${\mathbb R}$)  denotes  the group of fractional linear transformations of the upper half plane, and its double cover SL(2,${\mathbb R}$) the area-preserving linear mappings of the plane ${\mathbb R}^2$. 
SL(2,${\mathbb Z}$) $<$ SL(2,${\mathbb R}$) is the subgroup that preserves the integral lattice in ${\mathbb R}^2$, and PSL(2,${\mathbb Z}$) is the projection of SL(2,${\mathbb Z}$) to PSL(2,${\mathbb R}$).  Riemann surfaces (with known elementary exceptions)  are quotients of  the upper half plane  by torsion free discrete subgroups of PSL(2,${\mathbb R}$), and lifts of these subgroups to SL(2,${\mathbb R}$) correspond  precisely to spin structures on these Riemann surfaces.
In the universal setting, the mapping class group is the collection ${\rm PPSL}(2,{\mathbb Z})$
of all piecewise ${\rm PSL}(2,{\mathbb Z})$ homeomorphisms of the unit circle $S^1$ with finitely many breakpoints among the rational points in $S^1$.   Elements of ${\rm PPSL}(2,{\mathbb Z})$ are automatically 
${\mathcal C}^1$-  but never ${\mathcal C}^2$-smooth at minimal breakpoints.
Moreover, there is a triple incarnation of this group also as the Richard Thompson
group $T$ \cite{dudes} and the combinatorial Ptolemy group acting on tesselations with doe by flips from \cite{unicon}.
The new spin mapping class group P(SL(2,${\mathbb Z}$))
consists of all piecewise-constant maps $S^1\to{\rm SL}(2,{\mathbb Z})$ which projectivize to an element of ${\rm PPSL}(2,{\mathbb Z})$.
This group is defined combinatorially by its action on ${\mathcal Tess}^+$ and then identified with
the spin mapping class group in the proof of the next result.

\medskip

\begin{figure}[hbt] 
\centering   
\includegraphics[trim =0 430 0 90,width=1.\linewidth]{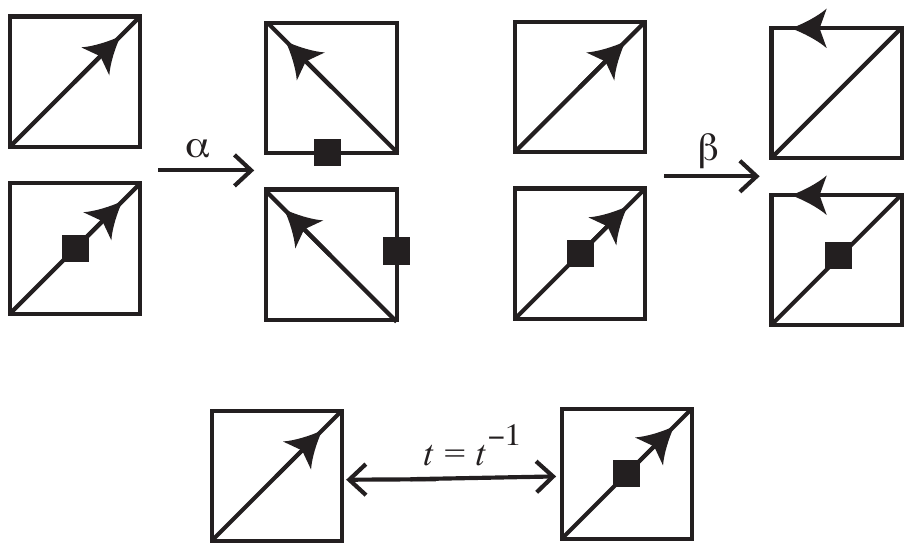}       
\caption{The three generators $\alpha,\beta, t$ of P(SL(2,${\mathbb Z}$)) acting on a marked tesselation with distinguished oriented edge 
in the two cases that the distinguished edge may or may not have a non-zero label in ${\mathbb Z}/2$ as indicated with a box icon. The distinguished oriented edge is indicated with an arrow.}
\label{fig:spin}               
\end{figure}

Here then is the birthday gift for Po:

\medskip

\noindent{\bf Main Theorem}.~{\it
The group ${\rm P(SL(2},{\mathbb Z}))$ generated by the
three transformations $\alpha,\beta,t$ illustrated in Figure \ref{fig:spin} 
acting on ${\mathcal Tess}^+$
is precisely the group of piecewise maps
$S^1\to{\rm SL(2,}{\mathbb Z}{\rm )}$ which projectivize to homeomorphisms
in {\rm PPSL(2,${\mathbb Z}$)}.}

\medskip

My goal here is to explain these new constructions 
with motivations going back to first principles in some detail, which I hope might be useful whether as survey or invitation.
The entire story comes together diagrammatically in Figure~\ref{fig:spin} illustrating the action of the generators of the universal spin mapping class group
on its Teichm\"uller space, a kind of emblem that I hope Milen might even appreciate as design.

A companion paper \cite{comp} to this one derives a finite presentation of 
P(SL(2,${\mathbb Z}$)) from that of PPSL(2,${\mathbb Z}$) $\approx T$.  In fact,
\cite{comp} contains the first complete derivation in the literature of the 
latter finite presentation as well, though an equivalent presentation was
given in \cite{LS} seemingly based on \cite{unicon} and on
unpublished notes of Thompson.  By construction,  the diagram
$$\begin{array} {c c c c}
{\rm P(SL(2},{\mathbb Z}))\times{\mathcal Tess}^+&\rightarrow&{\mathcal Tess}^+\\
{~~~~~~~~~~~~~~~~~}\downarrow&&\downarrow\\
{\rm PPSL}(2,{\mathbb Z})\times{\mathcal Tess}&\rightarrow&{\mathcal Tess}\\
\end{array}$$
commutes, where the vertical maps are the natural forgetful maps.

As I first learned
as a graduate student from the Orsay {\it Travaux de Thurston sur les Surfaces} seminars \cite{FLP} co-organized by Po, the classical Teichm\"uller theory originated as a topic in complex analysis and evolved
under William Thurston's masterful infusion of hyperbolic geometry.  
This has led 
to a combinatorial approach useful more recently for studying the Riemann moduli space and Teichm\"uller space $T(F)$ of an orientable surface
$F=F_g^s$ of genus $g\geq 0$ with $s\geq 1$ punctures, where $2-2g-s<0$.

The key construction, as recalled later, canonically associates a decomposition $\Delta(\Gamma)$ of $F$ into ideal polygons
to a hyperbolic 
structure $\Gamma$ on $F$ suitably decorated with one horocycle about each of its punctures.
Furthermore, {\sl solving the equation $\Delta(\Gamma)=\Delta$ in $\Gamma$ for fixed $\Delta$} yields an ideal cell decomposition of the decorated Teichm\"uller space itself, which is a difficult
theorem in \cite{penner}.  Generically, an ideal triangulation of $F$ is the label for a top-dimensional cell
in the decorated Teichm\"uller space of $F$.
In effect in the universal setting, an ideal triangulation of the surface is replaced by a tesselation  of its universal cover.

For a {\sl universal Teichm\"uller space}, we demand \cite{unicon} a Fr\'echet manifold ${\mathcal T}$ supporting
the action ${\mathcal P}\acts{\mathcal T}$ of a {\sl universal mapping class group} ${\mathcal P}$ so that:

\medskip

\leftskip .3in
\rightskip.3in

\noindent $\bullet$~there are injective homomorphisms of the classical mapping class groups $MC(F)$ into explicit completions $\widehat{\mathcal P}_F\acts{\mathcal T}$ of ${\mathcal P}$;

\smallskip

\noindent $\bullet$~there are embeddings $T(F)\subset{\mathcal T}$
of the classical Teichm\"uller spaces which are equivariant for the $MC(F)$-action;

\smallskip

\noindent $\bullet$~the Weil-Petersson K\"ahler geometry on each $T(F)$ is induced by pulling back
a universal geometry on ${\mathcal T}$.

\medskip
\leftskip=0ex\rightskip=0ex

Our universal Teichm\"uller space satisfying these conditions and more
is given by the
quotient $${\mathcal Tess}={\mathcal Tess}'/{\rm PSL}(2,{\mathbb R}).$$
Suitably decorating a tesselation with horocycles at its ideal points, a construction 
analogous to the classical case for finite-type surfaces (namely, a convex hull construction in Minkowski space,
cf.~the appendix) provides 
an ideal polygonal decomposition of ${\mathbb D}$, and the space of decorated tesselations
itself again admits a corresponding decomposition by {\sl solving the equation}.

Also from \cite{unicon}, there are global so-called
{\sl lambda length coordinates} on ${\mathcal Tess}$ providing the Fr\'echet structure.
In fact, the log lambda length deformations also from \cite{unicon}
lead to an infinite-dimensional Lie algebra of piecewise ${\frak sl}_2$
vector fields on $S^1$ with breakpoints among the rationals,
a kind of completion of the loop algebra of ${\frak sl}_2$, introduced in \cite{MP}
and more recently studied in \cite{FP} as part of a larger program with implications
to and from the current paper.
This Lie algebra  denoted  ${\rm p} {\frak sl}_2$ in \cite{MP} (for piecewise  ${\frak sl}_2$) was renamed
${\rm pp} {\frak sl}_2$ in \cite{FP} (to emphasize its relation with ${\rm PPSL}(2,{\mathbb Z})$) and
 is now more precisely seen to be {\sl exactly the Lie algebra of} P(SL(2,${\mathbb Z}$)).
 Lambda lengths and the associated Lie algebra are not further discussed in the current paper.

The novel contributions of this paper involve including spin structure in the universal setting.
But it is not so much novel as it is the
natural universal \cite{unicon} extension  from \cite{N=1} in the classical case.
The latter is reformulated in \cite{N=2},
as we shall explain, with related aspects of GL(1$|$1)-graph connections also studied in \cite{andrea}.
An overview of the several equivalent formulations of spin structures on finite-type surfaces is given in the appendix.
Considerations of universal spin structure are prefatory to any discussion of super universal Teichm\"uller theory.

The most elegant and immediately useful characterization of spin structures in our setting is due to Natanzon \cite{natanzon}
for
a hyperbolic structure on $F$ specified by a (conjugacy class of projective) uniformizing representation $\pi_1(F)\to{\rm PSL}(2,{\mathbb R})$
of the fundamental group: 
a {\it spin structure on $F$} is a
lift $\pi_1(F)\to{\rm SL}(2,{\mathbb R})$ of the uniformizing representation from
${\rm PSL}(2,{\mathbb R})$  to  ${\rm SL}(2,{\mathbb R})$. 

The natural spin generalization of ${\rm PPSL}(2,{\mathbb Z})$ based on Natanzon's description is our {universal spin mapping class group} P(SL(2,${\mathbb Z}$))
of all piecewise maps $S^1\to{\rm SL}(2,{\mathbb Z})$ with no conditions at the breakpoints 
except we demand that the underlying quotient map $S^1\to{\rm PSL}(2,{\mathbb Z})$ lies in
${\rm PPSL}(2,{\mathbb Z})$.  Composition in the group P(SL(2,${\mathbb Z}$)) is given by taking common refinements of the piecewise structures and ordinary composition in ${\rm SL}(2,{\mathbb Z})$ on the resulting pieces, just as in ${\rm PPSL}(2,{\mathbb Z})$.

To define the spin
universal space, we extend a different characterization of spin structure from \cite{N=1}  given by equivalence classes of orientations $\omega:\tau\to\tilde\tau$  on the edges in $\tau\in{\mathcal Tess}$, where
$\tilde\tau$ denotes the set of oriented edges of $\tau$. The equivalence relation is 
generated by a move called {\it (Kastelyn) reflection} on an ideal triangle complementary to $\tau$, which reverses the orientations of all three frontier edges of a fixed complementary triangle.  

In fact, the doe $e$ itself also extends to another orientation $\omega_e:\tau\to\tilde\tau$ induced by the flow in $S^1$ from initial to terminal point on its right and the reverse orientation on its left.  (Any specific orientation will suffice to specify the 
``standard'' spin structure determined by the doe.)
It follows that on a tesselation $\tau$ with doe $e\in\tilde\tau$, we might specify the orientation
$\omega$ determining a universal spin structure by instead marking with a box icon each edge $f\in\tau$
where $\omega(f)\neq\omega_e(f)$, called a {\it marking}.  In this notation taken modulo two on each edge, a reflection
on a complementary ideal triangle simply adds a mark to each frontier edge of the triangle on which the reflection is performed.
Equivalently as in \cite{N=2}, \cite{andrea} and the Abstract, a spin structure is a ${\mathbb Z}/2$-graph connection on the graph dual
to an ideal triangulation of a punctured surface in the classical case, and also now dual to a tesselation with doe in the universal case. 
In this paper, we shall primarily stick to the formalism of equivalence classes of ${\mathbb Z}/2$-markings on edges of a tesselation.

$\S$\ref{over1} gives an overview of 
Teichm\"uller theory, both classical and universal, necessary for 
a further discussion of our results, with more complete details
but few proofs in an appendix.
(Let us apologize for pedagogically inevitable minor expositional repetition in slightly different contexts in these two surveys.)
$\S$\ref{unisp} discusses universal spin structure and the proof of the Main Theorem.  The reference \cite{pb} subsumes the classical \cite{penner} and universal \cite{unicon} decorated Teichm\"uller theories but does not treat spin structures from \cite{N=1}.

\section{overview of ${{\rm Pt}\approx\rm PPSL}(2,{\mathbb Z})\acts{\mathcal Tess}\approx{\rm Homeo}_+(S^1)/{\rm PSL}(2,{\mathbb R})$}\label{over1}

Fix an oriented surface $F=F_g^s$ as in the Introduction and choose a point in its Teichm\"uller space $T(F)$, i.e., specify some complete finite-area metric on $F$ of constant Gauss curvature -1 modulo push-forward by diffeomorphisms of $F$ isotopic to the identity. The universal setting is intended to provide an infinite-dimensional space ${\mathcal T}$, together with embeddings
of each $T(F)\subset {\mathcal T}$ which may depend upon choices, so that both the geometry and mapping class quotient topology of the classical spaces $T(F)$ pull-back corresponding universal structures on ${\mathcal T}$.

\begin{figure}[hbt] 
\centering          
\includegraphics[width=0.6\textwidth]{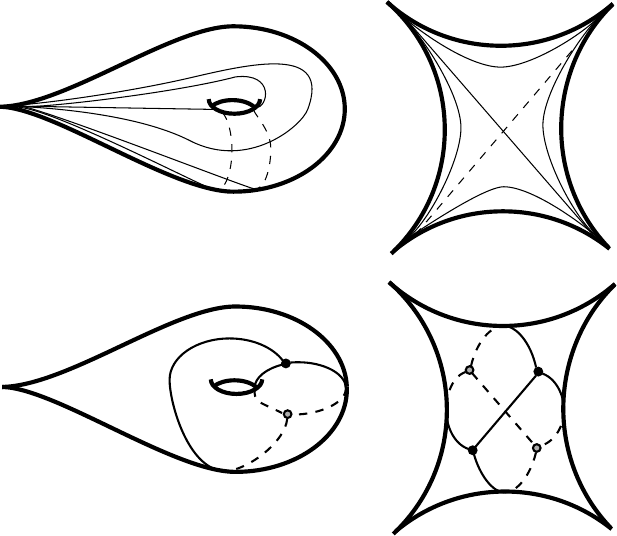}
\caption{Examples of ideal triangulations in $F_1^1$ and $F_0^4$ above and their dual fatgraphs below.}  
\label{fig:11}               
\end{figure}

It turns out that the further assignment of one real parameter to each puncture of $F$, which may be interpreted as
the hyperbolic length of a specified horocycle about the puncture, is enough to determine
a collection in $F$  of pairwise disjointly embedded arcs with endpoints among the punctures and decomposing $F$ into ideal polygons, as illustrated on the top in Figure~\ref{fig:11}.  (This extra decoration is also sufficient to define the global affine coordinates  on the decorated bundles given by  lambda lengths.) The existence of such
a canonical decomposition in Theorem \ref{icd} from \cite{penner} is the starting point for all the combinatorics.

The Poincar\'e dual of this cell decomposition in $F$ is a graph $G$ embedded as its deformation retract,
as illustrated on the bottom of Figure \ref{fig:11}, and it has the extra structure of a cyclic ordering on the half-edges about each vertex
coming from an orientation on $F$, a so-called {\it fatgraph} $G\subset F$.  Generically in Teichm\"uller space and as in the figure, the decomposition $\Delta$ is an {\it ideal triangulation} with only complementary triangles, and in this case $G$ is a trivalent fatgraph.

In effect, the passage to the universal case entails lifting such an ideal triangulation
of $F$ to 
its universal cover ${\mathbb D}$, namely the open unit disk in the complex plane 
${\mathbb C}$ with the Poincar\'e metric.  In the generic case, the ideal triangulation 
$\Delta$ of $F$ lifts to a {tesselation} $\tau$ of ${\mathbb D}$. Of course, such a lift to ${\mathbb D}$ is invariant under the action of the cofinite Fuchsian group uniformizing $F$.  The universal case amounts to keeping the tesselation but dropping the Fuchsian group.  

 \begin{figure}[hbt!] 
\centering   
\includegraphics[trim =0 60 0 35,width=.75\linewidth]{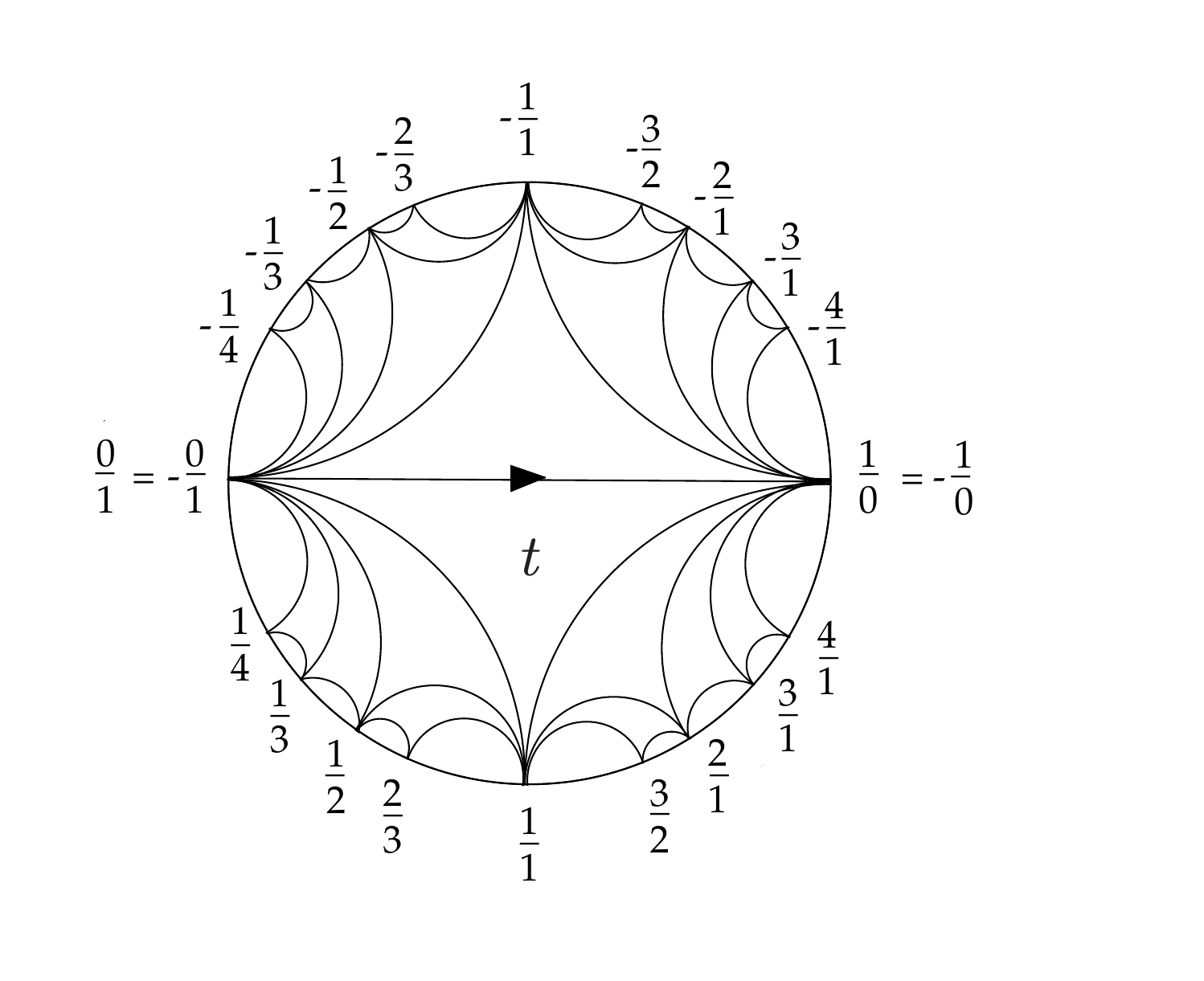}  
\caption{The Farey tesselation $\tau _*$ of the Poincar\'e disk ${\mathbb D}$ with its triangle $t$ to the right of its distinguished oriented edge.  The extended rationals ${\mathbb Q}\cup\{\infty\}$ 
are naturally identified with the ideal vertices of $\tau_*$, as
illustrated in low generation.}  \label{farey}               
\end{figure}

In order to have one natural example in mind and as discussed in detail in the appendix, the {\it Farey tesselation} $\tau_*$  is illustrated in Figure \ref{farey} with its {\it distinguished oriented edge} or {\it doe} $e_*\in\tilde\tau_*$,
where $\tilde\tau$ denotes the set of orientations on edges in 
$\tau$.  
Define $${\mathcal Tess}'=\{{\rm tesselations~with~doe~of~}{\mathbb D}\}.$$

The doe is not intrinsic and its specification as extra data is critical to the universal setting since
tesselations with doe are rigid in the following sense:  Let $\tau^0\subset S^1=\partial {\mathbb D}$ denote the collection of ideal points of edges in $\tau$.  Given two tesselations ${\tau}_0$ and $\tau$ with doe,  there is a unique $f:\tau_0^0\to\tau^0$ so that $x,y\in S^1$ are ideal vertices of an edge in $\tau_0$ if and only if $f(x),f(y)$ are likewise for $\tau$.  This is easy to see inductively,
one triangle at a time starting with the triangle to the right of the doe.  As explained in the appendix, this map $f$ moreover interpolates an orientation-preserving homeomorphism also denoted $f:S^1\to S^1$.

Taking the Farey tesselation $\tau_0=\tau_*\in{\mathcal Tess}'$ with doe  as a chosen base point in this construction gives the {\it characteristic map}
$$f_{e\in\tilde\tau}\in{\rm Homeo}_+(S^1)=\{{\rm orientation{-}preserving~homeomorphisms~of~}S^1\}$$ canonically associated to $e\in\tilde\tau$, and in fact, this association of characteristic map $f_{e\in\tilde\tau}$
to a tesselation $\tau$ with doe $e$ induces an isomorphism
$${\mathcal Tess}'\approx{\rm Homeo}_+(S^1).$$

\begin{figure}[hbt] 
\centering          
\includegraphics[width=0.55\textwidth]{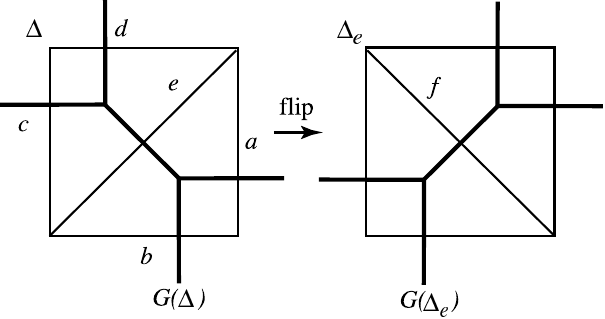}
\caption{The flip on an edge $e$ in ideal triangulation $\Delta$ produces another
ideal triangulation $\Delta_e$, and equivalently the flip in its dual trivalent fatgraph 
$G(\Delta)$ produces $G(\Delta_e)$.}  \label{flip0}               
\end{figure}

Our model of {\it universal Teichm\"uller space} is the quotient
$${\mathcal Tess}={\mathcal Tess}'/{\rm PSL}(2,{\mathbb R})\approx{\rm Homeo}_+(S^1)/{\rm PSL}(2,{\mathbb R}),
$$
and its {\it universal mapping class group} is generated by a combinatorial move, called a {\it flip}
defined as follows for any edge $e\in\tau\in{\mathcal Tess}$:
remove $e$ from $\tau$ so as to produce a complementary quadrilateral
and then replace the one diagonal $e$ with the other diagonal $f$ of this quadrilateral, as indicated in Figure~\ref{flip0}, along with the dual move on fatgraphs.  Using the characteristic map,
edges and hence flips can be indexed by elements $\tau_*$, so the collection of all finite sequences\footnote{
The completion 
$\widehat{\mathcal P}_F$ of PPSL(2,${\mathbb Z}$) in the wish-list for a universal theory in the Introduction arises from finite sequences of simultaneous flips along all the edges in an orbit of the uniformizing Fuchsian group, i.e., the lift to the universal cover
of a flip on one edge in the underlying surface.}
of
such elements forms a group, the {\it Ptolemy group} Pt.

In fact via the characteristic map, the Ptolemy group Pt can be identified with
the subgroup ${\rm PPSL}(2,{\mathbb Z})$
of ${\rm Homeo}_+(S^1)$ consisting of
all piecewise-${\rm PSL}(2,{\mathbb Z})$ homeomorphisms of $S^1$ with only 
finitely many pieces whose endpoints lie in $\tau_*^0$, and this group turns out to
furthermore be isomorphic to Richard Thompson's group $T$, all of which is detailed in the appendix. Thus, our universal mapping class group
$${\rm Thompson~group}~T \approx {\rm PPSL}(2,{\mathbb Z})
\approx {\rm Ptolemy~group~Pt}$$
acts on the universal Teichm\"uller space ${\mathcal Tess}$ by flips.

\section{Proof of Main Theorem}
\label{unisp}

One key explanation for the gift/theorem/emblem in the Introduction is that the {\sl enhanced flip on an oriented edge in a tesselation}
engendered by the move $\alpha$ in Figure \ref{fig:spin} is precisely the dual of the flip transformation on fatgraphs with spin structure
discovered and illustrated in Theorem \ref{t:N=1} from \cite{N=1,N=2}.

We have already defined the universal spin mapping class group 
$${\rm P(SL(2},{\mathbb Z}))=\Biggl\{
\begin{aligned}
&{\rm piecewise}~\phi:S^1\to{\rm SL}(2,{\mathbb Z})~{\rm with~rational~breakpoints}\\
&{\rm so~that~the~projectivization~of~}\phi{\rm ~lies~in}~{\rm PPSL}(2,{\mathbb Z})
&
\end{aligned}
 \Biggr\}
 $$
and universal spin Teichm\"uller space ${\mathcal Tess}^+={\mathcal Tess'}^+/{\rm PSL}(2,{\mathbb R})$, where
 $${\mathcal Tess'}^+=\{(e\in\tilde\tau,\mu)\in{\mathcal Tess}'\times\{ 0,1\}^\tau:\mu~{\rm has~finite~support}\}/\sim$$
with the equivalence relation $\sim$ on {markings} $\mu$ given by finite compositions of  {\it (Kastelyn) reflections} which
change marking on all three frontier edges of any fixed complementary triangle.  
Just as in \cite{unicon} for flips, infinite sequences of reflections can also be allowed provided the marking of each edge eventually stabilizes.

For $A=\begin{psmallmatrix}a&b\\c&d\\\end{psmallmatrix}\in{\rm SL}(2,{\mathbb R})$, define the sign
$${\rm sgn}(A)=\begin{cases}
{\rm sgn}(a+d);&a+d\neq 0,\\
{\rm sgn}(c);&{\rm else},
\end{cases}$$
where 
${\rm sgn}(x)$ denotes the usual sign of $0\neq x\in{\mathbb R}$.
In particular,  ${\rm trace} (A)=a+d=0=c$ is impossible since $ad-bc=1$.
It follows that a spin hyperbolic structure on a finite-type surface may be thought of as a pair comprised of a projective uniformizing representation
$$\pi_1(F)\to{\rm PSL}(2,{\mathbb R})$$ as usual together with 
an induced function $${\rm sgn}:\pi_1(F)\to\{\pm 1\}$$ determining which sign
to assign to the representing projective matrix, the one of positive trace (for sgn=1) or of negative trace (for sgn=-1)
generically, with a similar definition using the sign of the (2,1)-entry of the representing matrix in the traceless case.

Analogously, given  $\phi:S^1\to{\rm SL}(2,{\mathbb Z})$ in P(SL(2,${\mathbb Z}))$,
there is by definition the induced projectivization $\bar\phi:S^1\to{\rm PSL}(2,{\mathbb Z})$
in ${\rm PPSL}(2,{\mathbb Z})$.
Using the sign function
in the previous paragraph, $\phi$ itself is uniquely determined by $\bar\phi$ together with a function
$\sigma:\Pi(\bar\phi)\to\{\pm 1\}$ determining the signs, where $\Pi(\bar\phi)=\Pi(\phi)$ is the set of pieces of $\bar\phi$ or $\phi$.


 \begin{figure}[hbt] 
\centering   
\includegraphics[trim =0 340 140 170,width=.60\linewidth]{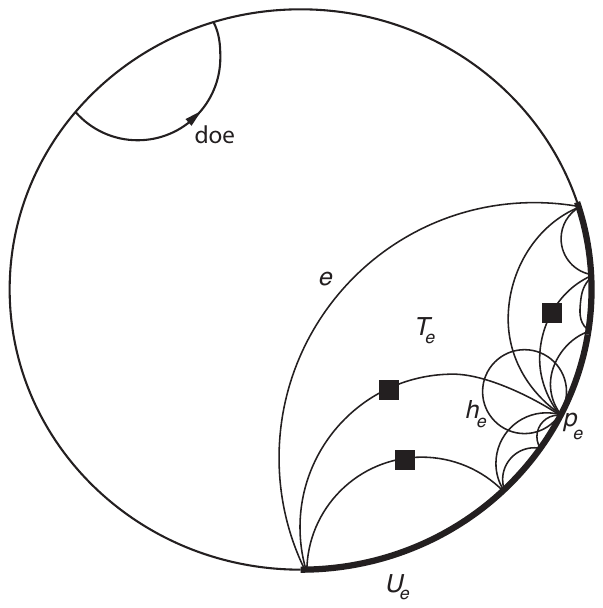}  
\caption{In the presence of the doe, any other edge $e\in\tau$ independently of orientations determines
a complementary triangle $T_e$ 
with its ideal vertex $p_e=\partial T_e-\partial e\in S^1$ in a component $U_e$ of $S^1-\partial e$ in boldface.
Also illustrated is a horocycle $h_e$ centered at $p_e$.  A finite number of nearby edges
may have non-zero markings, as indicated with the box icon.
}  \label{phoro}               
\end{figure}

\begin{figure}[tbhp]
\centering
\includegraphics[trim = 0 420 0 260,width=1.25\linewidth]{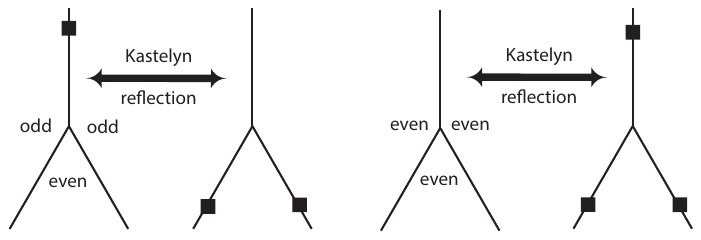}
\caption{ 
Reflection preserves parity modulo two of the number of
marked edges on each nearby horocycle as depicted on the level of dual fatgraphs.}
\label{fig:K}
\end{figure}

To define $\sigma$, we first define an auxiliary function
$\sigma:\tau-\{doe\}\to\{\pm 1\}$ as follows.   As depicted in Figure \ref{phoro} in the presence of a doe, another edge
$e\in\tau$ determines a point $p_e\in U_e\subset S^1$, where $U_e$ is the component of $S^1-\partial e$ on the other side from the doe.  A
horocycle $h_e$ 
centered at $p_e$ meets a countable collection of edges of $\tau$, some finite number of which have non-zero marking,
and we define $\hat\sigma(e)=+1$ if this number is even, and otherwise $\hat\sigma(e)=-1$, completing the definition
of $\hat\sigma$.  As illustrated in Figure \ref{fig:K}, this function $\hat\sigma:\tau\to\{\pm1\}$ is well-defined on
reflection equivalence classes.  Let us emphasize that this determination according to parity of the number of marked
edges meeting a horocycle derives from the classical case of finite-type surfaces.  It is not simply invented here for our separate
purposes.  

Now suppose that $e\in\tilde\tau$, where $\tau$ agrees with the Farey tesselation $\tau_*$ outside of a finite
ideal polygon $P$ which moreover contains all edges with non-zero markings as well as the doe.  An induction over innermost disks
proves that there is a
marking equivalent to the given one with all non-zero marked edges lying in the frontier of $P$, so we may henceforth assume
that the marking is of this type.

The ideal vertices of $P$
contain the breakpoints of $\Pi(\phi)$, and it therefore  suffices to define $\sigma$ on the pieces so determined, and we
set $\sigma (I)=\hat\sigma(e)$ if the circular interval $I$ comprising a piece  in the piecewise structure of $\phi$ agrees with some $U_e$.
This defines the mapping $$\langle \alpha,\beta,t\rangle\to{\rm P(SL(2,}{\mathbb Z}))$$ from the group generated by $\alpha,\beta,t$ to the group of piecewise SL(2,${\mathbb Z}$)-valued maps which descend to PPSL(2,${\mathbb Z}$).

It remains to prove that $\langle \alpha,\beta,t\rangle\to$~P(SL(2,${\mathbb Z}$)) is bijective.  We begin with the
proof of surjectivity, and
to this end claim that there is a finite composition of reflections on triangles in $P$ achieving any specified
marking on its frontier, which we take to be an ideal $n$-gon.  To see this, consider the fatgraph $G$ dual to $P$, 
which has $V_e=n$ univalent external vertices and $V_i=n-2$ trivalent internal ones, with
$n$ edges incident on the former and $n-3$ edges with endpoints among the latter, for a total of $E=2n-3$,
so of course $V_e+V_i-E=1$
according to Euler characteristic, whence $V_e=1+E-V_i$.

Rather than Kastelyn reflect on a complementary ideal triangle in $P$, let us dually imagine reflecting on an interior  
vertex of $G$ by changing orientations on each incident edge, and consider
the collection of all functions from the set of all edges of $G$ to ${\mathbb Z}/2$, a set with cardinality $2^{E}$. 
Interior vertex reflections act on this set of functions in the natural way, and there are evidently $2^{V_i}$ possible such compositions. 
The simultaneous reflection at all interior vertices of $G$ acts trivially on this set of functions, so only
$2^{V_i-1}$ such compositions may act non-trivially.  
Insofar as $2^E/2^{V_i-1}=2^{1+E-V_i}$ and there are
$1+E-V_i$ external edges, the claim follows, proving surjectivity.

As for injectivity of $\langle \alpha,\beta,t\rangle\to$~P(SL(2,${\mathbb Z}$)), consider an element of the
kernel, namely a word $w$ in $\alpha,\beta,t$ which in particular descends to the identity
in PPSL(2,${\mathbb Z}$) upon setting $t=1\in$~P(SL(2,${\mathbb Z}$)). In other words
in the language of the companion paper \cite{comp}, $w$ is a $t$-insertion in a relation
in PPSL(2,${\mathbb Z}$)).  Since it is in the kernel, $w$ moreover preserves some marking, and hence
any marking since markings, just as spin structures, provide a torsor.  This is precisely the characterization
of a relation in P(SL(2,${\mathbb Z}$)) based on orbifold covers of quotients that is exploited in \cite{comp} to derive the finite presentation
of P(SL(2,${\mathbb Z}$)).
An element in the kernel of $\langle \alpha,\beta,t\rangle\to$~P(SL(2,${\mathbb Z}$)) is thus a relation
in the domain, proving injectivity.
~~~\hfill $\qed$

\appendix

\section{background material}

\subsection{Farey tesselation and modular group}\label{modgrp}

Let ${\mathbb Z}\subseteq{\mathbb Q}\subseteq{\mathbb R}\subseteq{\mathbb C}$ denote the integers, rational, real and complex numbers, respectively, 
and let $i=\sqrt{-1}$.  Let $\infty$ indicate a {\it point at infinity}
and the superscript plus sign denote
the one-point completion by $\infty$, e.g., ${\mathbb Q}^+={\mathbb Q}\cup\{\infty\}$.
The upper half-plane ${\mathcal U}=\{
x+iy\in {\mathbb C}:y> 0\}$ endowed with the metric  $ds^2=\frac{dx^2+dy^2}{y^2}$ provides
a standard model for hyperbolic geometry.  
The {\it Cayley transform} 
$C:
s\mapsto {{s-i}\over{s+i}}$ induces an isomorphism of pairs $({\mathcal U},{\mathbb R}^+)\to ({\mathbb D},S^1)$,
where
${\mathbb D}\approx{\mathcal U}$ is the open unit disk in the complex plane with the
induced metric, or {\it Poincar\'e disk} model of the hyperbolic plane, with its frontier $S^1$ the {\it circle at infinity}.

Let $t$ denote the ideal
triangle with vertices
$\pm1, -i\in S^1$ as in Figure~1, and consider the
group $\Gamma$ generated by hyperbolic reflections in its frontier edges.  Define
the {\it Farey tesselation} $\tau _*$ to be the $\Gamma$-orbit
of the frontier edges of $t$, regarded as a set of geodesic {\it edges},  
together with its {\it distinguished
oriented edge}, or {\it doe}, given by the edge from
$-1$ to $+1$. 
 A rational point ${p\over q}\in {{\mathbb Q}^+=\mathbb Q}\cup\{\infty\}$, as illustrated in Figure~1, is said to be of {\it generation g} if the radial arc in ${\mathbb D}$ from the origin to ${p\over q}$ meets the interior of $g\geq 0$ distinct complementary triangles.

The {\it modular group} ${\rm PSL}(2,{\mathbb Z})$
is the subgroup of $\Gamma$
consisting of compositions of an even number of reflections, or in other words
the group of two-by-two integral matrices $A$ of unit determinant modulo
the equivalence relation generated by identifying $A$ with $-A$.  This group acts on the left
by orientation-preserving hyperbolic isometry on $z\in{\mathcal U}$
by fractional linear transformation
$A=\begin{psmallmatrix} a&b\\c&d\\\end{psmallmatrix}:z\mapsto {{az+b}\over{cz+d}}$.  There is also
the action of $A\in {\rm PSL}(2,{\mathbb Z})$  {on the right} (following Gauss) on the rational points,
for which we introduce the notation
$A:{\rm doe}\mapsto {e}_A\in\tilde\tau$, where
$e_A$ has respective initial and terminal points ${{b+ia}\over{b-ia}}$ and ${{d+ic}\over{d-ic}}$ in ${\mathbb C}$, that is,
respective labels $-{b\over a}$ and $-{d\over c}$ in Figure~\ref{farey}.

The main point for us is that the modular group leaves set-wise invariant the Farey tesselation $\tau_*$, mapping $\cup\tau_*$ onto $\cup\tau_*$, and any orientation-preserving homeomorphism of the circle leaving $\tau_*$ invariant in this manner lies in ${\rm PSL}(2,{\mathbb Z})$. In fact, the modular group acts simply transitively on the set
$\tilde\tau_*$ of orientations on edges of $\tau_*$, and
a generating set for ${\rm PSL}(2,{\mathbb Z})$ is given by any two of 
$$R=\begin{pmatrix}0&-1\\1&~1\\\end{pmatrix}, ~
S=\begin{pmatrix}0&-1\\1&~0\\\end{pmatrix}, ~
T=\begin{pmatrix}1&-1\\0&~1\\\end{pmatrix}, ~
U=\begin{pmatrix}1&0\\1&1\\\end{pmatrix}. $$
Moreover,  $S^2 = 1 = R^3$ is a complete set of relations in the generators
$R=T^{-1}U$ and $S=TU^{-1}T$, 
so ${\rm PSL}(2,{\mathbb Z}) \approx {\mathbb Z}/2\star{\mathbb Z}/3$. 

\begin{figure}[hbt] 
\centering          
\includegraphics[width=0.30\textwidth]{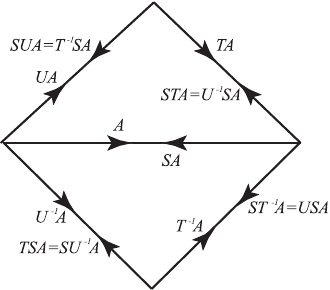}
\caption{Enumeration of oriented edges near $e_A$ in $\tau_*$.}  \label{label}               
\end{figure}

Geometrically in Figure~1: the elliptic element
$S$ fixes the distinguished edge in $\tau_*$ reversing the orientation of the doe;
$R$ is the elliptic transformation which counter-clockwise cyclically permutes the vertices
of the triangle $t$; and $U$ (respectively  $T$) is 
the parabolic transformation
with the fixed point ${0\over 1}$ (respectively ${1\over 0}$)
which cyclically permutes the incident edges of $\tau _*$ in
the counter-clockwise sense about ${0\over 1}$ (respectively
the clockwise sense about ${1\over 0}$).  Typical aspects of this enumeration
of oriented edges by elements of ${\rm PSL}(2,{\mathbb Z})$ are illustrated in Figure~\ref{label}. 

\subsection{Tesselations of ${\mathbb D}$}
An arbitrary tesselation of ${\mathbb D}$ is a collection $\tau$ of bi-infinite geodesics 
locally finite in ${\mathbb D}$ and decomposing it into complementary ideal triangles. 
Geodesics in $\tau$ are called its {\it edges}, and $\tau$ itself is regarded as a set of edges. Let $\tau^0\subset S^1$ denote the collection of endpoints of all the edges in $\tau$, which is automatically dense in $S^1$ by local finiteness using that complementary regions are ideal triangles.
Let $\tilde\tau$ denote the set of oriented edges of $\tau$.  
A distinguished oriented edge, or simply {\it doe}, on $\tau$ is the specification of an element of $\tilde\tau$. 

Another canonical tesselation of ${\mathbb D}$ with doe is the dyadic tesselation $\tau_d$ , which has the same doe as $\tau_*$ as well as the same generation zero vertices, and which is recursively characterized by the property that one vertex of each triangle complementary to $\tau_d$ bisects the angle between its other two vertices. Thus, $\tau_d^0$ consists of the points in $S^1$ with (dyadic) rational arguments, as opposed to
$\tau_*^0\subset S^1$ which is comprised of points with rational rectilinear coordinates.

We again exploit the fact that tesselations with doe are rigid, namely, suppose that $\tau$ is any tesselation of 
${\mathbb D}$ with doe $e\in\tilde\tau$.  There is a unique $A\in{\rm PSL}(2,{\mathbb R})$ 
mapping the doe of $\tau_*$ to $e$ and the triangle $t$ to its right to the triangle complementary to $\tau$ on the right of $e$.  One continues in this manner mapping complementary triangles to induce a map
$\phi=\phi_{e\in\tilde\tau}:\tau_*^0\to\tau^0$ called the {\it characteristic map} of $e\in\tilde\tau$.

It is not hard to see that $\phi$ is a surjection using local finiteness in ${\mathbb D}$, and it is an order-preserving injection by construction.  As an order-preserving bijection between dense sets in $S^1$, the characteristic map interpolates an orientation-preserving homeomorphism $\phi:S^1\to S^1$ of the same name.  Letting 
${\rm Homeo}_+(S^1)$ denote the topological group of orientation-preserving homeomorphisms of $S^1$ with the compact-open topology and ${\mathcal Tess}'$ the space of tesselations of ${\mathbb D}$ with doe with the Hausdorff topology on closed subsets of ${\mathbb D}$, we have

\begin{theorem}[\cite{unicon}] \label{homeo}The assignment $\phi_{\tau,e}$ of characteristic map
to $e\in\tilde\tau$ induces a homeomorphism ${\mathcal Tess}'\to{\rm Homeo}_+(S^1)$.
\end{theorem}

\noindent Our model of {\it universal Teichm\"uller space} is the quotient
$${\mathcal Tess}={\mathcal Tess}'/{\rm PSL}(2,{\mathbb R})\approx{\rm Homeo}_+(S^1)/{\rm PSL}(2,{\mathbb R})
,$$
which can be identified with the collection of all ideal triangulations of ${\mathbb D}$
with the same doe as $\tau_*$ and the same triangle $t$ to its right; in effect in ${\mathcal Tess}'$,
the doe determines a triangle, which then is normalized in ${\mathcal Tess}$ to $t$ by ${\rm PSL}(2,{\mathbb R})$.
Bers' universal Teichm\"uller space \cite{bers-univ} is the collection of all quasi-symmetric homeomorphisms of $S^1$ modulo ${\rm PSL}(2,{\mathbb R})$, so our construction naturally generalizes Bers' version by extending to all orientation-preserving homeomorphisms of $S^1$.

\subsection{Teichm\"uller theory for finite-type surfaces}

Consider a connected closed orientable surface with genus $g\geq 0$, and let 
$F=F_g^s$ be the complement  of a finite set of cardinality $s\geq1$, which are taken to be the punctures of $F$.  Let us tacitly fix an orientation on $F$ and choose a basepoint for its fundamental group.

Provided the Euler characteristic $2-2g-s<0$ is negative,  
$F={\mathcal U}/\Gamma$ is {\it uniformized by a Fuchsian group}, by which we mean: let 
${\mathcal U}$ denote the upper half plane
with its Poincar\'e metric
 $ds^2=\frac{dx^2+dy^2}{y^2}$ and its projective matrix group ${\rm PSL}(2,{\mathbb R})={\rm SL}(2,{\mathbb R})/\pm I$ of oriented isometries, where $I$ denotes the identity matrix;
upon choosing a basepoint in $F$, there is a (conjugacy class of) injective homomorphism called the {\it uniformizing representation} $\pi_1\to {\rm PSL}(2,{\mathbb R})$ of the fundamental group $\pi_1=\pi_1(F)$
of $F$ onto a discrete {\it Fuchsian subgroup} $\Gamma<{\rm PSL}(2,{\mathbb R})$ so that non-trivial loops about punctures in $F$ are represented by parabolic transformations, namely, those with absolute trace equal to two.
See  \cite{beardon,ford,pb} for example.

The {\it Teichm\"uller space} of $F$ is
$$T(F)~=~{\rm Hom}'(\pi_1,{\rm PSL}(2,{\mathbb R}))/{\rm PSL}(2,{\mathbb R}),$$
where the prime indicates Fuchsian representations as just defined and the action
of ${\rm PSL}(2,{\mathbb R})$ on ${\rm Hom}'$ is by conjugation.  
The {\it decorated Teichm\"uller space} is simply $\tilde T(F)=T(F_g^s)\times{\mathbb R}_+^s$, where the
decoration is conveniently regarded as an $s$-tuple of positive real weights on the punctures.
In particular, the {\it mapping class group} $MC(F)$ of homotopy classes 
of orientation-preserving homeomorphisms of $F$ acts on $T(F)$ and on $\tilde T(F)$ in the natural way by push-forward.
$T(F)$ is homeomorphic to an open ball of real dimension $6g-6+2s$.

The homotopy class of a collection $\Delta$ of pairwise disjointly embedded essential arcs in $F$ connecting punctures is called an {\it arc family} in $F$ provided no two distinct arcs in $\Delta$ are isotopic.  $\Delta$ is said to {\it fill} $F$ provided each component of $F-\Delta$ is simply connected.  An {\it ideal cell decomposition}
is a decomposition of a topological space into simplices minus certain of their faces of codimension at least two, e.g., an ideal triangulation $\Delta$ filling $F$ decomposes it into triangles with their {\sl ideal vertices at infinity} at the punctures of $F$.

There is the following foundational result based on a convex hull construction\footnote{
Uniformize the punctured surface $F$ in Minkowski space by $\Gamma\acts{\mathbb R}^{2,1}$, where $\Gamma$
acts by Lorentz isometry on ${\mathbb R}^{2,1}$ and lies in the component of the identity.  A choice of horocycles about
the punctures  in $F$ provides a finite collection of $\Gamma$-orbits of points in the positive light-cone in
${\mathbb R}^{2,1}$, the closed convex hull of which provides a $\Gamma$-invariant convex body in ${\mathbb R}^3\approx {\mathbb R}^{2,1}$.
The extreme edges of this body project to an ideal cell decomposition of $F$.  This basic convex hull construction is from \cite{penner}.
The same applies to the analogous end for discrete radially dense subsets of the light-cone in the universal context of \cite{unicon}.}  in Minkowski space:

\begin{theorem}[\cite{penner,pb}] \label{icd} There is an $MC(F)$-invariant smooth ideal cell decomposition ${\mathcal C}(F)$ of
$\tilde T(F)$  whose simplices are indexed by arc families filling $F$
with faces corresponding to inclusion of arc famiies.
\end{theorem}

A maximal arc family is called an {\it ideal triangulation} of $F=F_g^s$ and contains $6g-g+3s$ edges.
In fact, the lambda lengths discussed in the Introduction on the edges of any ideal triangulation give
global affine coordinates on $\tilde T(F)$.
By construction, crossing a codimension-one face of ${\mathcal C}(F)$ corresponds to a {\it flip} on an edge $e$ in an ideal triangulation $\Delta$ as illustrated in Figure \ref{flip0}, replacing $e$ by $f$ to produce another ideal triangulaton $\Delta_e$ of $F$.

As also depicted, dual to an arc family $\Delta$ filling $F$, there is a 
graph $G=G(\Delta)\subset F$ embedded in $F$ as a deformation retract, called a {\it spine} of $F$.
An orientation on $F$ induces the counter clockwise ordering on the half edges of $G$ incident on each fixed vertex thus giving the abstract graph $G$ the structure of a 
{\it fatgraph} (also called a {\it ribbon graph}). 
Dual to flipping diagonals of a quadrilateral in an ideal triangulation of $F$, there is the
combinatorial move depicted in Figure \ref{flip0} also called a {\it flip} on the
dual trivalent fatgraph spine $G\subset F$:
contract an edge of $T$ with distinct endpoints and then expand the resulting 4-valent vertex in the unique distinct manner in order to produce another trivalent fatgraph.   

This leads to the {\it Ptolemy groupoid} ${\rm Pt}(F)$ of $F$ whose objects are homotopy classes of  ideal triangulations of $F$, or dually trivalent fatgraph spines, and whose morphisms are compositions of flips. ${\rm Pt}(F)$ is the fundamental path groupoid of $\tilde T(F)$ according to Theorem \ref{icd}.
In fact, finite compositions of flips act transitively on homotopy classes of trivalent fatgraph spines.
To see this, use the evident path connectivity of $\tilde T(F)$
to join by a smooth path cells corresponding to distinct ideal triangulations. Put this path in general position
with respect to the ideal cell decomposition ${\mathcal C}(F)$ so that it crosses only codimension-one faces, as required.
See \cite{penner,pb} for the details, namely, ${\mathcal C}(F)$ is algebraic hence smooth so that general position is viable,
and every ideal triangulation of $F$ actually occurs in ${\mathcal C}(F)$.

It follows that flips generate $MC(F)$ in the sense that if $G\subset F$ is a trivalent fatgraph spine and $\varphi\in MC(F)$, then there is a sequence
$$\varphi(G)=G_1-G_2-\cdots -G_n=G$$ of trivalent fatgraph spines of $F$ where any
consecutive pair differ by a flip, with the similar statement for ideal triangulations.   Moreover, the general position argument just given, but now for 
two-dimensional homotopies of paths, implies

\begin{figure}[hbt] 
\centering           
\includegraphics[width=0.9\textwidth]{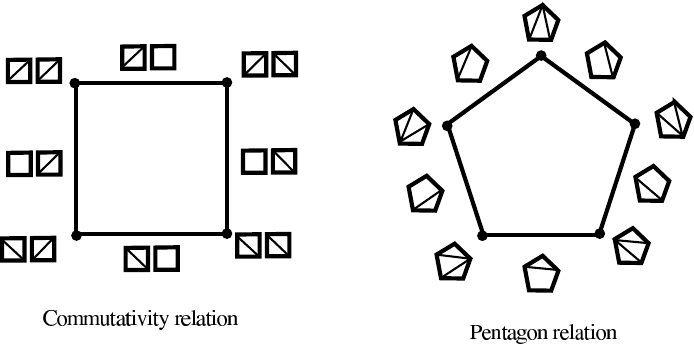}
\caption{Links of codimension-two cells in the idea cell decomposition of decorated Teichm\"uller spaces.} \label{fig:31}               
\end{figure}

\begin{corollary}\label{relations}
A complete set of relations in the fundamental path groupoid of $\tilde T(F)$ is given by 
declaring flips idempotent together with
the links of codimension-two faces as illustrated in Figure~\ref{fig:31}, namely: commutativity of flips
on disjoint quadrilaterals; and the pentagon relation of five consecutive flips alternating between two distinct edges
of a triangle.
\end{corollary}


\subsection{Ptolemy group}
The {\it Ptolemy group(oid)} Pt has objects given by tesselations with doe of ${\mathbb D}$ which coincide with $\tau_*$ outside of a finite ideal polygon and morphisms given by finite compositions of flips.  
A typical flip clearly has order two, however, the flip on the doe is defined to produce 
another doe 
so that these oriented edges in this order respect the orientation of $F$, and hence the flip on the doe has order four, cf.~the top-left for $\alpha$ in Figure \ref{fig:spin} ignoring boxes

Using combinatorial rigidity of
tesselations with doe, flips can be labeled by edges of a fixed tesselation, so words in these labels render Pt in fact a group, not just a groupoid  as it is for finite-type surfaces $F$ without doe.   Furthermore, ${\rm PSL}(2,{\mathbb Z})$ sits inside Pt as those tesselations which are identical to $\tau_*$ except perhaps for the location of the doe.

Let ${\rm PPSL}(2,{\mathbb Z})<{\rm Homeo}_+(S^1)$ denote the collection of all piecewise-${\rm PSL}(2,{\mathbb Z})$ homeomorphisms of $S^1$ with only finitely many pieces whose endpoints lie in $\tau_*^0$.
Since any orientation-preserving homeomophism of ${\mathbb D}$ preserving $\cup\tau_*$ must lie in the modular group, it follows that ${\rm PPSL}(2,{\mathbb Z})$ can also be described as the collection of characteristic
maps of objects in Pt.
As pointed out to me by Maxim Kontsevich \cite{Ko} decades ago now and as exposed by Imbert in \cite{imbert}, 
by the same logic these two discrete subgroups of ${\rm Homeo}_+(S^1)$ coincide under the isomorphism in Theorem \ref{homeo}, so we have ${\rm PPSL}(2,{\mathbb Z})\approx{\rm Pt}$.

Moreover, it is not difficult to see that the characteristic map $\tau_*\to\tau_d$ of the dyadic tesselation $\tau_d$ is precisely the Minkowski ?-function, and it
conjugates ${\rm PPSL}(2,{\mathbb Z})$ to the celebrated Thompson group $T$, cf.~\cite{dudes}.  We have arrived at the remarkable fact that the three avatars 
$${\rm Thompson~group}~T \approx {\rm PPSL}(2,{\mathbb Z})
\approx {\rm Ptolemy~group~Pt}$$
of our {\it universal mapping class group} Pt act on the universal Teichm\"uller space ${\mathcal Tess}\approx
{\rm Homeo}_+(S^1)/{\rm PSL}(2,{\mathbb R})$ by flips.  

Let $\alpha$ denote the flip on the doe and
$\beta$ denote the move that fixes the tesselation and
moves the doe around the triangle to its left, cf.~$\beta$ in Figure \ref{fig:spin}
again ignoring boxes.
Clearly $\alpha^2\sim S$ and $\beta\sim R$ generate
${\rm PSL}(2,{\mathbb Z})$, which
acts simply transitively on $\tilde\tau_*$.

As explained in \cite{comp}, there is the following direct consequence of the ideal analogue of Theorem \ref{icd} from \cite{unicon}:

\begin{theorem} {\rm PPSL(2,${\mathbb Z})$} is generated by the flip $\alpha$ on the doe and the transformation $\beta$ which moves the doe one edge counter-clockwise in the triangle to its left.  A presentation in these generators is given by the following relations:
$\alpha^4$, $\beta^3$, $(\alpha\beta)^5$ and the two commutators $[\beta\alpha\beta,\alpha^2 \beta\alpha\beta \alpha^2]$ and
$[\beta\alpha\beta,\alpha^2 \beta^2\alpha^2~\beta\alpha\beta~\alpha^2\beta\alpha^2]$.
\end{theorem}

\noindent See \cite{LS} for the discussion of an algebraic proof of an equivalent presentation based on unpublished computations of Thompson, and \cite{comp} for a complete and self-contained proof.

\subsection{Spin structures on finite-type surfaces} 
\label{over2}

Milnor's elegant definition \cite{milnor} of a {\it spin structure} on a smooth surface $F$ is a class in the modulo two first cohomology group 
of the unit tangent bundle 
of $F$ which is non-zero on the fiber class.
An equally elegant definition of immediate utility in our situation due to Natanzon \cite{natanzon} on  a uniformized surface $F$
is given by a lift of the uniformizing representation $\pi_1\to {\rm PSL}(2,{\mathbb R})$ to $\pi_1\to {\rm SL}(2,{\mathbb R})$;
this immediately leads to P(SL(2,${\mathbb R}$)) as the universal spin mapping class group.

\begin{figure}[hbt] 
\centering           
\includegraphics[width=0.65\textwidth]{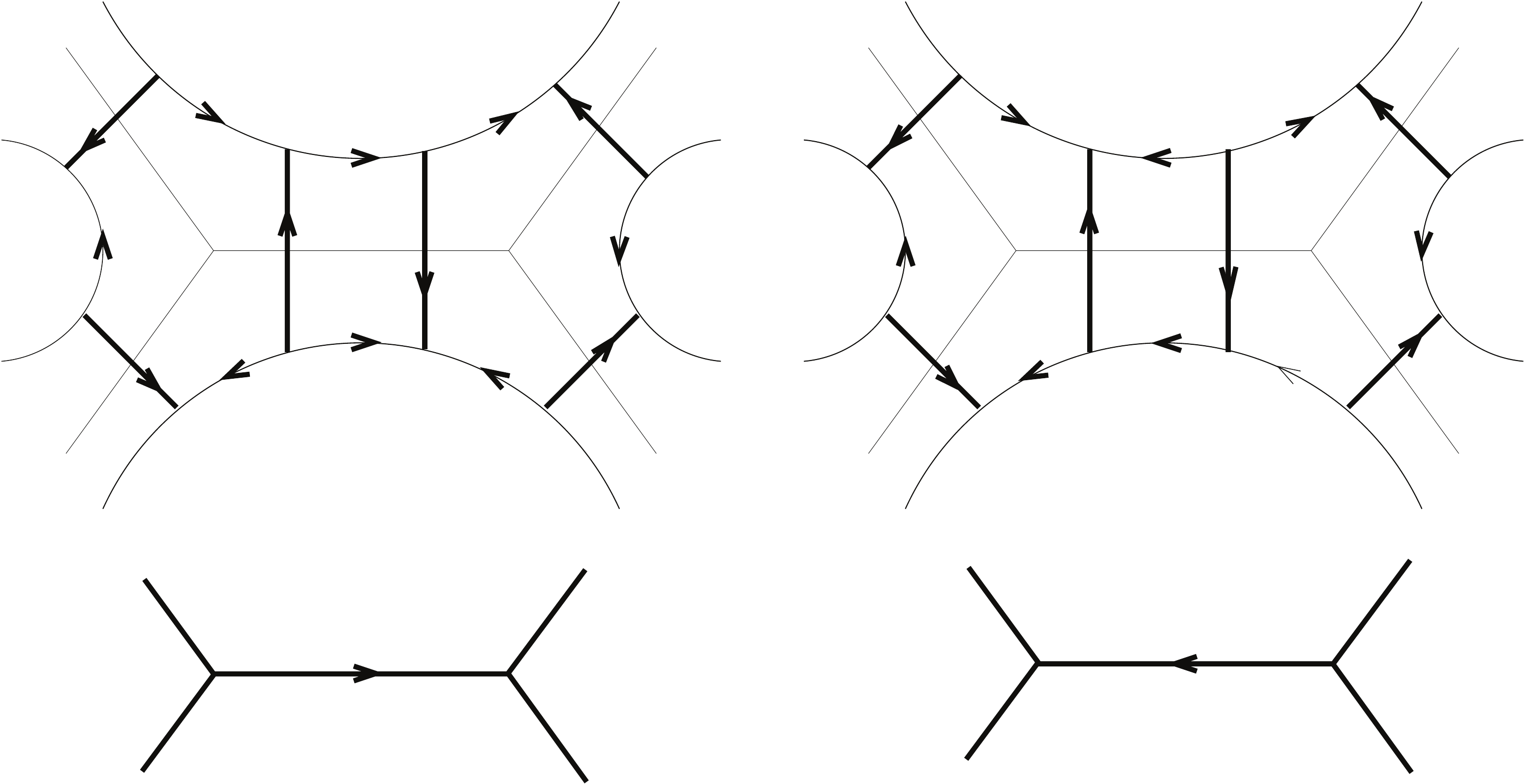}
\caption{On the top, CW decomposition of a neighborhood of a trivalent fatgraph $G$ with its special dimer indicated in boldface.  Orientations on hexagons extend to Kastelyn orientations on the top, and on the bottom is the corresponding  orientation on the edge of $G$ in the two cases.} \label{K}               
\end{figure}

Our starting point for the universal spin Techm\"uller space  is Johnson's general formulation  \cite{johnson}
as an element of the affine $H^1(F;{\mathbb Z}/2)$-space ${\mathcal Q}(F)$
 of quadratic forms,
i.e., functions $$q:H_1=H_1(F;{\mathbb Z}/2)\to{\mathbb Z}_2$$
which are quadratic with respect to the homology intersection pairing 
$\cdot:H_1\otimes H_1\to {\mathbb Z}_2$ on
$H_1$ in the sense that
$q(a+b)=q(a)+q(b)+a\cdot b$, for $a,b\in H_1$.
Our combinatorial starting point towards Johnson's definition as an element of ${\mathcal Q}(F)$
is the special case of a general technique 
 from statistical physics
by Cimazoni-Reshetihkin in \cite{cr1,cr2}.  This is expressed in terms of Kastelyn orientations (which disagree with the
orientation as a boundary an odd number of times, cf.~the top of Figure~\ref{K}) and dimer configurations (disjoint unions of closed edges covering the vertices, again cf.~Figure \ref{K})
on the one-skeleton of a suitable CW decomposition of $F$:

\begin{theorem}[Theorem 2.2 of \cite{cr2}]\label{cr}
Fix a dimer configuration $D$ on a surface graph with boundary $\mathcal{G}$ for  the surface $\Sigma$ and let 
$\alpha\in H_1(\Sigma; \mathbb{Z}_2)$ be represented by oriented closed curves $C_1,\dots, C_m\in \bar{\mathcal{G}}$. If $K$ is a Kasteleyn orientation on $\mathcal{G}$, then the function $q^K_D: H_1(\Sigma; \mathbb{Z}_2)\to \mathbb{Z}_2$ given by 
\begin{eqnarray}\label{q}\nonumber
q^K_{D}(\alpha)=\sum_{i<j}C_i\cdot C_j+\sum^m_{n=1}(1+n^K_{C_i}+\ell^D_{C_i})\quad  ({\rm mod} ~2)
\end{eqnarray}  
is a well-defined quadratic form,
where $\ell^D_C$ is the number of edges of $D$ sticking out to the left of $C$,
and
$n^K_C$ is the number of edges counted with multiplicity where the orientation of $C$ disagrees with that of $K$.  Moreover for each fixed dimer $D$, this establishes
an isomorphism as affine $H^1(\Sigma; \mathbb{Z}_2)$-space between 
${\mathcal Q}(\Sigma)$ and the collection of equivalence classes of Kastelyn orientations
on $\mathcal{G}$, with  equivalence generated by reversing orientations
of all edges incident on some fixed vertex of ${\mathcal G}$.
\end{theorem}

We actually rely here on another formulation of
spin structure from \cite{N=1} which amounts to choosing a CW decomposition in Theorem \ref{cr} suited to the tesselation with
a special dimer on it, cf. Figure \ref{K}.  
Fix a trivalent fatgraph spine $G\subset F$
and consider the set ${\mathcal O}(G)$ of equivalence classes of all orientations
on its edges, where the equivalence relation is 
generated by {\it reflections} at vertices $v$ of by reversing the orientation of each edge incident on $v$.
(These reflections on $G$ are compositions of six Kastelyn reflections
on $\overline{\mathcal G}$ as in Theorem \ref{cr}.)

Here is the combinatorial result from \cite{N=1} which explains the genesis of the generator 
$\alpha$ of P(SL(2,${\mathbb Z}$)):


\begin{theorem}[\cite{N=1,N=2}]\label{t:N=1}
Suppose $G\subset F=F^s_g$ is a trivalent fatgraph spine.  Then $\mathcal{O}(G)$ and  $\mathcal{Q}(F)$ are isomorphic as affine $H^1(F; \mathbb{Z}_2)$-spaces. Moreover, the action of the mapping class group $MC(F)$ on ${\mathcal Q}(F)$
lifts to the action of the Ptolemy groupoid on $\mathcal{O}(G)$ described by the following figure
\medskip

\begin{tikzpicture}[ultra thick, baseline=1cm,scale=.7]
\draw (0,0)--(210:1) node[above] at (210:0.7){$\e_2$};
\draw (0,0)--(330:1) node[above] at (330:0.7){$\e_4$};
    \draw[ 
 	ultra thick,
        decoration={markings, mark=at position 0.5 with {\arrow{>}}},
        postaction={decorate}
        ]
        (0,0) -- (0,2);
\draw[yshift=2cm] (0,0)--(30:1) node[below] at (30:0.7){$\e_3$};
\draw[yshift=2cm] (0,0)--(150:1) node[below] at (150:0.7){$\e_1$};
\end{tikzpicture}
\begin{tikzpicture}[baseline]
\draw[->, thick](0,0)--(2,0);
\node[above] at (1,0) {};
\node at (-1,0){};
\node at (3,0){};
\end{tikzpicture}
\begin{tikzpicture}[ultra thick, baseline,scale=.7]
\draw (0,0)--(120:1) node[above] at (100:0.3){$\e_1$};
\draw (0,0)--(240:1) node[below] at (260:0.3){$\e_2$};
    \draw[ 
        decoration={markings, mark=at position 0.5 with {\arrow{<}}},
        postaction={decorate}
        ]
        (0,0) -- (2,0);
\draw[xshift=2cm] (0,0)--(60:1) node[above] at (100:0.3){$-\e_3\hskip1ex$};
\draw [xshift=2cm](0,0)--(-60:1) node[below] at (-80:0.3){$\e_4$};
\end{tikzpicture},

\medskip
\noindent where $\epsilon_i$ indicates orientations on edges, and the extra minus sign denotes orientation reversal;
in the special case $F=F_1^1$, $\e_3$ replaces $-\e_3$.
\end{theorem}

\noindent Embarrassingly, there were inadvertently two cases after a lengthy reduction from 16 cases in \cite{N=1}, but as noted in \cite{N=2}, the two cases are each reflection-equivalent to the single diagram in the theorem.

A notion from physics, a spin structure  $q\in{\mathcal Q}(F)$ on $F$ distinguishes two types of punctures: If $\gamma_p$ is a simple loop about the puncture $p$ of $F$, then $p$ is a {\it Neveu-Schwarz} puncture if $q ([\gamma_p])=0$ and otherwise $p$ is a {\it Ramond} puncture, where
bracket denotes homology class.

\begin{theorem}[\cite{ramond}] Consider a fatgraph spine $G$ of $F$ and a simple oriented edge-path $\gamma$ in $G$ surrounding a fixed puncture $p$ of $F$.  A spin structure as in the previous theorem determined by the  equivalence class of an orientation $\omega$ on $G$ has $p$ as a Ramond puncture if and only if the length of  $\gamma$ has the same parity modulo two as the number of edges of $G$ where $\omega$ agrees with $\gamma$. \end{theorem}

\vfill\eject


\begin{thebibliography}{ABCD}



\bibitem{beardon} Beardon, A., {\it The Geometry of Discrete Groups}, Springer-Verlag, Berlin 1983.

\bibitem{bers-univ}
Bers, L.,
 Universal Teichm\"uller space,
In {\it Analytic Methods in Mathematical Physics}, Gordon and Breach,
New York, 1968, 65-83.

\bibitem{Po++}
Boone, W.W., Haken, W.,  Po\'enaru, V.,
On Recursively Unsolvable Problems in Topology and Their Classification
{\it Studies in Logic and the Foundations of Mathematics},
volume 50, 1968, Pages 37-74,  Elsevier Publishers, Amsterdam Netherlands.

\bibitem{andrea}
Bourque, A. and Zeitlin, A.,
Flat GL(1$|$1)-connections and fatgraphs
preprint arXiv:2208.08033 (2022)

\bibitem{dudes}
Cannon, J., Floyd, W., Parry, W., Introductory remarks on Richard Thompson’s groups. L’Enseignement Mathematique 42, 215–256 (1996)


\bibitem{cr1}
Cimasoni, D. and  Reshetikhin, N.,
 Dimers on surface graphs and spin structures.\ I, {\it Communications in Mathematical Physics} {\bf 275} (2007), 187-208.

\bibitem{cr2}---,  Dimers on surface graphs and spin structures.\ II, {\it Communications in Mathematical Physics} {\bf 281} (2008), 445-468.

\bibitem{FLP}
Fathi, A., Laudenbach, F., Po\'enaru, V.,
{\it Travaux de Thurston sur les Surfaces} Asterisque {\bf 66-67},
Soc. Math. de France, Paris, 1979.


\bibitem{ford} Ford, L.,
{\it  Automorphic Functions}, American Mathematical Society (1951).

\bibitem{FP}
Frenkel, I. and Penner, R.,
 Sketch of a Program for Automorphic Functions from Universal Teichm\"uller Theory to Capture Monstrous Moonshine, {\it Communications in Mathematical Physics} {\bf 389} (2022), 1525-1567.
 
 \bibitem
{Po+}
Haefliger, A. and  Po\'enaru, V.,
 La classification des immersions combinatoires, 
 {\it Publications Math\'ematiques de l'IH\'ES} {\bf 23} (1964), 75–91.


\bibitem{imbert}
Imbert, M.,
Sur l'isomorphisme du groupe de Richard Thompson avec le groupe de Ptolme
 {\it Geometric 
Galois Actions II}, London Math Society Lecture Notes {\bf 243}, 
Cambridge University Press (1997), eds.\ P.\ Lochak and L.\ Schneps.


\bibitem {N=2} Ip,I., Penner, R., Zeitlin, A., {N}=2 Super-Teichm\"uller theory,{\it Advances in Mathematics} {\bf 336} (2018), 409-454.

\bibitem{ramond} ---,  On Ramond Decorations,  {\it Communications in Mathematical Physics} {\bf 371} (2019), 145-157.


\bibitem{johnson} Johnson, D.,
Spin structures and quadratic forms on surfaces, {\it J.\ London Math.\ Soc.}{\bf 22} (1980) 365-373. 

\bibitem
{Ko} Kontsevich, M.,
Private communication 1995.


\bibitem{LS}
Lochak, P., Schneps, L.,
The universal Ptolemy-Teichm\"uller groupoid,
 {\it Geometric 
Galois Actions II}, London Math Society Lecture Notes {\bf 243}, 
Cambridge University Press (1997), eds.\ P.\ Lochak and L.\ Schneps.

\bibitem{MP}
Malikov, F. and Penner, R.,
The Lie algebra of homeomorphisms of the
circle, {\it Advances in Mathematics} {\bf 140} (1999), 282-322.

\bibitem{barry}
B. Mazur, A note on some contractible 4-manifolds, {\it Annals of Math.} {\bf 73} (1961),
221-228.




\bibitem{milnor} Milnor, J.,
 Remarks concerning spin manifolds in {\it Differential and Combinatorial Topology}, Princeton University Press, 1965.

\bibitem{natpaper} Natanzon, S.M.,
 Moduli of Riemann surfaces, Hurwitz-type spaces, and their superanalogues, {\it Russian Math.\ Surveys} {\bf 54} 1 (1999) 61-117.

\bibitem{natanzon}---, 
Moduli of Riemann Surfaces, Real Algebraic Curves, and Their Superanalogs, {\it Transl.\ Math.\ Monographs} {\bf 225},
American Mathematical Society, 2004.




\bibitem{penner}Penner, R.C.,
The decorated Teichm\"uller space of punctured surfaces,
{\it Communications in Mathematical Physics} {\bf 113} (1987), 299-339.


\bibitem{unicon}
---, Universal constructions in Teichm\"uller theory,  
{\it Advances in Mathematics}  {\bf   98} (1993), 143-215.


\bibitem{pb} ---, {\it Decorated Teichm\"uller theory}, European Mathematical Society, 2012. 

\bibitem{comp} ---,
Universal Spin Teichm\"uller Theory, II. Finte Presentation of P(SL(2,${\mathbb Z}$)), preprint, 2023.

\bibitem{N=1}
Penner,  R. and Zeitlin, A.,
Decorated super-Teichm\"uller space,
{\it Journal of Differential Geometry} {\bf 111} (2019), 527-566.

\bibitem{QSF1}
Po\'enaru, V.,
Geometric simple connectivity and low-dimentional topology
(2003)
mathnet.ru/links/f1937b66240c8471d3a2533073eed2f8/tm19.pdf


\bibitem{QSF2}
---,
Geometric simple connectivity and finitely presented groups
arxiv.org/abs/1404.4283
(2014)

 \bibitem{QSF3}
---,
 All finitely presented groups are QSF
arxiv.org/abs/1409.7325
 (2014)
 
  \bibitem{3+4}
 ---,
 On the 3-Dimensional Poincaré Conjecture and the 4-Dimensional Smooth Schoenflies Problem, 
arxiv.org/abs/math/0612554
(2006)


\bibitem{pocell}
---, La decomposition de l'hypercube en produit topologique, {\it Bull. Soc.
Math. France} {\bf 88} (1960), 113-129.




 



\end{thebibliography}
\end{document}